\documentclass[11pt]{amsart}

\title{On linearity of automorphism groups of algebras and groups  }

\author{Vitali\u\i\ Roman'kov}
\address{Institute of Mathematics and Information Technologies\\Dostoevsky Omsk State  University}
\curraddr{}
\email{romankov48@mail.ru}

\usepackage{amsmath,amsthm,amsfonts,amssymb}
\usepackage{graphicx}
\usepackage{hyperref}
\usepackage[usenames]{color}

\newtheorem{theorem}{Theorem}[section]

\theoremstyle{definition}

\newtheorem{corollary}[theorem]{Corollary}

\newtheorem{problem}[theorem]{Problem}
\newtheorem{remark}[theorem]{Remark}

\newcounter{comcount}

\date{}

\begin{document}

\maketitle

\begin{abstract}
The main purpose of this paper is to describe some published results and outline corresponding approaches which when applied to automorphism groups of algebras or groups establish that these groups are linear or non-linear.
\end{abstract}

\section{ On linearity of automorphism groups of groups}
\label{se:intro}

There are some useful tests for linearity, such as
\begin{enumerate}
\item  A linear group has ascending chain condition on centralizers.
\item (Mal’cev) Finitely generated  linear groups are residually finite.
\item  (Mal’cev) A solvable linear group is nilpotent-by-(abelian-by-finite).
\item (Tits) A linear group either contains a free group of rank two, or is solvable-by-(locally finite).
\end{enumerate}
\subsection{Non-linear automorphism groups of finitely generated solvable groups}
In \cite{T}, Tits showed that a finitely generated matrix group either contains a solvable normal subgroup of a finite index (i.e., is almost solvable) or else contains a noncyclic free subgroup.
In \cite{BM}, Bachmuth and Mochizuki conjectured that Tits' alternative is satisfied in every finitely generated group of automorphisms of a finitely generated solvable group. They point out that their conjecture holds for solvable abelian-by-nilpotent  groups and in some other cases. It turned out that this conjecture breaks down, in general. It has been independently shown by Hartley \cite{H} and the author \cite{Rom1} with absolutely different approaches.
We will give these results and corresponding approaches starting with \cite{Rom1}.
\begin{theorem}
 (Roman’kov \cite{Rom1}). The direct wreath product of the group of IA-automorphisms of an arbitrary finitely generated solvable group with an infinite cyclic group is embeddable into the group of automorphisms of some finitely generated solvable group.
\end{theorem}

Recall that the IA-automorphisms of a group $G$ are those automorphisms which are identical modulo the commutator group $G'$.
\begin{corollary} If $A$ is an arbitrary finitely  generated  almost solvable group, then the wreath product $B =A wr C_{\infty}$ ($C_{\infty}$ denotes the infinite cyclic group) is embeddable into the group of automorphisms of some finitely generated  solvable group.
\end{corollary}
The deduced corollary yields a negative answer to the Bachmuth and Mochizuki's question since, under its conditions, the Tits' alternative is satisfied in the group $G$ if and only if $A$ is solvable. Note that $A$  can be chosen as arbitrary finite group. 

In \cite{BM}, Bachmuth and Mochizuki also conjectured that in every finitely generated  group of automorphisms of a finitely generated solvable group, there is  a subnormal series of a finite length whose factors are either abelian or  matrix groups over a commutative Noetherian ring.
\begin{corollary}
\label{co:2}
 There exists a finitely generated group of automorphisms of a finitely generated solvable group which does not have the indicated property.
\end{corollary}

The infinite countable direct power $G^{\infty}$ of a simple finite non-abelian group $G$ is not embeddable into a group with the indicated subnormal series. It is easy to prove this by induction on the length of the series using the remark of Merzlyakov that the direct product of an infinite number of non-abelian groups cannot be represented by matrices (see \cite{M}, Example 3). The proofs for fields and for commutative Noetherian rings are identical. This power $G^{\infty}$ is the basic subgroup of $G wr C_{\infty}$. By Corollary \ref{co:2}, the latter group is embeddable in the group of automorphisms of some finitely generated solvable group.

Recall, that a group is called {\it perfect} if it equals its  derived subgroup. In \cite{H}, a group $G$ is called {\it perfectly distributed}, if every subgroup of a finite index of $G$ contains a non-trivial finitely generated perfect subgroup. Clearly no perfectly distributed group can be soluble-by-finite.
\begin{theorem}
\label{th:2}
 (Hartley \cite{H}). There exists a finitely generated solvable group $G$ of derived length three whose automorphism group contains subgroups $K \leq L$ such that $L$ is finitely generated, $L/K$ is infinite cyclic, $K$ is perfectly distributed and locally finite.

Moreover, it can be even arranged that $K$ is locally a direct power of any finite non-abelian simple group. Clearly $L$ is not solvable-by-finite, nor does it contain a non-abelian free subgroup.
\end{theorem}
The group $G$ of Theorem \ref{th:2} is actually an extension of a locally finite group by an infinite cyclic group. But a somewhat more complicated version of the construction for Theorem \ref{th:2} gives
\begin{theorem}
\label{th:3}
3 (Hartley \cite{H}). There exists a finitely generated  solvable group $G$ of derived length four whose automorphism group contains a torsion-free subgroup $L$ having a normal subgroup $K$ such that $L$ is finitely generated, $L/K$ is infinite cyclic, $K$ is perfectly distributed, every finitely generated subgroup of $K$ is abelian-by-finite.
\end{theorem}

\subsection{On linearity of automorphism groups of finitely generated free groups}
Let $F_n$ be an absolutely free group with basis $X_n = \{x_1, ..., x_n\}.$ 
\begin{theorem}

(a) (Krammer \cite{K}) Aut $F_2$ is linear.

(b) (Formanek and Procesi \cite{FP}). Aut $F_n$ is not linear for $n \geq 3$.
\end{theorem}
The proof of (b) in \cite{FP} uses the representation theory of algebraic groups to show that a kind of diophantine equation between the irreducible representations of a group $G$ is impossible unless $G$ is abelian-by-finite. This leads to statement, which says that the HNN-extension
\begin{equation}
\label{eq:1}
H(G) = < G \times G, t| t(g,g)t−1 = (1,g) \ \textrm{for all} \   g \in G>
\end{equation}
\noindent
cannot be a linear group if $G$ is not nilpotent-by-abelian-by-finite. The statement (b) then is proved by showing that for $n ≥ 3$, the automorphism group of a free group $F_n$  contains $H(F_2)$.

Refer as a {\it unitriangular} automorphism of $F_n$ with respect to $X_n$ to every automorphism $\varphi$ determined by a mapping of the form
\begin{equation}
\label{eq:}
\varphi  : x_1\mapsto x_1, x_i \mapsto v_ix_i \  \textrm{for} \  i = 2, …, n, 
\end{equation}
\noindent
where $v_i = v_i(x_1, ..., x_{i-1}$  is an arbitrary element of $F_{i-1}$ = gp$(x_1, ..., x_{i-1})$. Every collection $(v_2 ,..., v_n) \in F_1 \times ... \times F_{n-1}$
 determines an automorphism $\varphi$ of $F_n$. The unitriangular automorphisms of $F_n$ constitute a subgroup UTAut$F_n$, which we call the {\it group of unitriangular automorphisms} of $F_n$. To simplify expressions, we denote it by $U_n$. It is easy to see that, up to isomorphism, $U_n$ is independent of the choice of $X_n$.
\begin{theorem} (Roman’kov \cite{Rom2}). The group $U_n$ of unitriangular automorphisms of the free group $F_n$  of rank  is linear if and only if $n \leq 3$.
\end{theorem}
As $U_1$ is trivial and $U_2$ is infinite cyclic, these groups are obviously linear. In \cite{Rom2}, it is shown that $U_3$ is isomorphic to a subgroup of Aut$F_2$. Then by Krammer's theorem \cite{K} $U_3$ is linear.

We need  two statements about any group $G$ that follow.
\begin{theorem}
(a) (Formanek, Procesi \cite{FP}). Let $\rho$ be a linear representation of $H(G)$. Then the image of $G \times  \{1\}$ has a subgroup of finite index with nilpotent derived subgroup, i.e., is nilpotent-by-abelian-by-finite.

(b) (Brendle, Hamidi-Tehrani \cite{BH-T}). Let $N$ be a normal subgroup of $H(G)$ such that the image of
$G \times \{1\}$ in $H(G)/N$ is not nilpotent-by-abelian-by-finite. Then $H(G)/N$ is not linear.
\end{theorem}
In [8], a subgroup $G$ of $U_4$ has been constructed that is isomorphic to a quotient $H(F_2)/N$ such that the image of $G \times \{1\}$ in $H(F_2)/N$ is not nilpotent-by-abelian-by-finite.

\subsection{On linearity of automorphism groups of finitely generated relatively free groups}
Given an arbitrary variety ${\mathcal C}$ of groups, denote by $G_n$ the free group in ${\mathcal C}$ with a fixed basis $Y_n = \{y_1, ..., y_n\}$. Refer as a {\it unitriangular} automorphism of $G_n$  with respect to $Y_n$ to every automorphism $\phi$  determined by a mapping of the form
\begin{equation}
\label{eq:3}
\phi : y_1\mapsto y_1, y_i\mapsto u_iy_i\  \textrm{for} \  i = 2, …, n, \end{equation}
\noindent where $u_i = u_i(y_1, ..., y_{i-1})$  is an arbitrary element of $G_{i-1} =$ gp$(y_1, ..., y_{i-1}).$  Every collection $(u_2,..., u_n) \in  G_1\times ... \times G_1 \times ... \times G_{n-1}$  determines an automorphism $\phi $  of $G_n$. The unitriangular automorphisms of $G_n$ constitute a subgroup UTAut $G_n$, which we call the {\it group of unitriangular automorphisms} of $G_n$. To simplify expressions, we denote it by $V_n$. It is easy to see that, up to isomorphism, $V_n$ is independent of the choice of $Y_n$.
\begin{theorem}
(Auslander and Baumslag \cite{AB}, Merzlyakov \cite{M}). The automorphism group of every finitely generated nilpotent-by-finite group is linear.
Moreover, the holomorph of a group of this type (hence, its automorphism group as well) admits a faithful matrix representation over the ring $\mathbb{Z}$  of integers.
\end{theorem}
 This implies in particular that the automorphism groups of a finite rank relatively free group $G_n$  in nilpotent-by-finite varieties ${\mathcal C}$ of groups admit faithful matrix representations.
\begin{theorem}
\label{th:8}
 (Olshanskii \cite{O}). If a relatively free group $G_n$ is neither free nor nilpotent-by-finite then Aut$G_n$ is not linear.
\end{theorem}
Observe also that Aut $G_{\infty}$  for every nontrivial variety ${\mathcal C}$ includes a subgroup consisting of all permutations of an infinite set $Y_{\infty}$  of free generators, which is isomorphic to the symmetric group on the infinite set. It is well known that the latter is not linear. Olshanskii’s proof of Theorem \ref{th:8} uses the following reasoning.

For every group $G$ there is a homomorphism $\sigma :  G \rightarrow$ Aut $G$ associating to each $h\in G$ the inner automorphism $\sigma (h): f \mapsto hfh^{-1}$. The kernel of $\sigma$ is the center C$(G)$ of $G$ and the image of $\sigma$  is the group Inn $G$ of inner automorphisms of $G$. The latter is normal in Aut $G$, which means that we may regard the quotient $G/{\textrm C}(G) \simeq {\textrm Inn} G$ as a normal subgroup of Aut $G$. 
It is shown in \cite{O}, that in the case of neither free nor nilpotent-by-finite relatively free group $G_n$ there exists an automorphism $\theta $ of $G_n$ such that the extension $P$ of $G_n/{\textrm C}(G_n)$ by means of $\theta $ is not linear. Furthermore, $\theta$ was chosen so that the eigenvalues of the induced linear transformation $(G_n)_{ab} = G_n/G_n'$ of abelianization, which in this case, on assuming that Aut $G_n$ is a linear group, is a rank $n$ free abelian group, were not roots of unity. The automorphism $\theta$ is certainly not unitriangular since all  eigenvalues of each unitriangular automorphism are equal to $1$. The group Inn $G_n/{\textrm C}(G_n)$ also fails to consist of unitriangular automorphisms. Thus, the result of Olshanskii and his method of proof give no information on the possible linearity of the group $V_n=$ UTAut$G_n$.

\begin{theorem}
(Erofeev and Roman’kov \cite{ER}).
\label{th:9}

(a) Let $G_n$ be a rank $n \geq  2$ relatively free group in an arbitrary variety ${\mathcal C}$ of groups. The following hold: $V_1$ is trivial, while $V_2$ is the cyclic group of the order equal to the exponent of ${\mathcal C}$. These groups admit faithful matrix representations. For $n \geq 3$, if $G_{n−1}$ is  nilpotent group then so is $V_n$. If $G_{n−1}$ is nilpotent-by-finite then $V_n$ is linear over $\mathbb{Z}$.

(b) Let $G_n$ be a rank $n \geq  3$ relatively free group in an arbitrary nontrivial variety ${\mathcal C}$ of groups distinct from the variety of all groups. If $G_{n−1}$ is not a nilpotent-by-finite group then the group $V_n =$ UTAut $V_n$ of unitriangular automorphisms of $G_n$ admits no faithful matrix representation over any field.
\end{theorem}
Thus, the claims of Theorem \ref{th:9} yield exhaustive information on the linearity of groups of unitriangular automorphisms of relatively free groups of finite rank in the proper varieties of groups. Namely, $V_n$ admits a faithful matrix representation over some field if and only if $G_{n−1}$ is a nilpotent-by-finite group.

\subsection{On linearity of automorphism groups determined by test elements}

 Let $G$ be a group. An element $g \in G$ is called a {\it test} element (for recognizing automorphisms) if, whenever $\varphi (g) = g$ for an endomorphism $\varphi $ of $G$, it follows that $\varphi $ is an automorphism. This definition was introduced by V. Shpilrain in \cite{Sh}. Let $T(G)$ denote the set of all test elements of $G.$

We say that a subgroup $H=H(g)$ of Aut$G$ is determined by the test element $g\in t(G)$ if  $H = \{\varphi \in  {\textrm{Aut}}G: \varphi (g) = g\}.$ In fact, $H(g)$ is a stabilizer of $g$ in Aut$ G.$ A question arises: what subgroups of Aut$G$ are determined by test elements? 

If we define a relation $\leq $ on the set $T(G)$ by $g\leq f \Leftrightarrow H(g) \subseteq H(f), $ a natural question arises about structure of $T(G)$ with the relation $\leq .$

More generally, we say that a subset $T \subset G$ is a {\it test} set if, whenever $\varphi (g) = g$ for every $g\in T$ for  an endomorphism $\varphi $ of $G$, it follows that $\varphi $ is an automorphism. The {\it test } rank $t(G)$ of $G$ is defined as the minimal size of a test subset of $G.$ Similarly, we introduce the relation $\leq $ on the set of all test subsets of $G.$

Now we just indicate some results about linearity of subgroups determined by test elements.  One of these results is implied by the following theorem. 

\begin{theorem} (Timoshenko \cite{Tim1}).
\label{th:1.12}
 Let $M_n$ be the free metabelian group of rank $n\geq 2$. Then

(a)  $T(M_2) =  M_2'\setminus \{1\}$; 

(b)  $t(M_n) = n-1$;

(c) a set of elements $\{g_1, ..., g_{n-1}\}$ is a (minimal) test set of $M_n$ if and only if, when  
$\{g_1, ..., g_{n-1}\}\subseteq M_n',$ they are linear independent over the group ring 
$\mathbb{Z}G/G'$ and every endomorphism $\phi $ fixing each of these elements induces an automorphism of $G/G'$. 
\end{theorem} 
\begin{corollary}
Every automorphism $\phi \in$ Aut$M_n$ determined by a minimal test set $\{g_1, ..., g_{n-1}\}$ mentioned fixes every element $u \in M_n',$ thus by Shmel'kin's theorem \cite{Shm} is inner corresponding to some element $h \in M_n'.$ 
\end{corollary}
It follows, that every subgroup $H=H(\{g_1, ..., g_{n-1}\})$ determined by a minimal test set $\{g_1, ..., g_{n-1}\}$ coincides with the subgroup Inn$_{M_n'}(M_n)$ of Inn$M_n$ corresponding to conjugations by elements of $M_n'.$ Hence, $H \simeq M_n'$ and clearly it  is a linear abelian group. 

It was shown in \cite{Romtest}, that in  the case of the free solvable group $S_{2,3}$ of rank $2$ and class $3$  there is a test element $g$ that lies in the second term $S_{2,3}^{(2)}$ of $S_{2,3}.$ Since this term has dimension $1$ over the group ring 
$\mathbb{Z}S_{2,3}/S_{2,3}'$ (see Remark below for explanation), it follows that the subgroup $H=H(g)$ of   Aut$S_{2,3}$ is isomorphic to $S_{2,3}^{(2)}$, hence is linear. By \cite{Tim2}, a similar assertions are true for any free solvable group $S_{2,n}$ of rank $2$.

\begin{remark}
Let $S_{n,d}$ be a free solvable group of rank $n\geq 2$ and class $d \geq 2$. Then the last non-trivial member $S_{r,d}^{(d-1)}$ can be treated as a module over the group ring $\mathbb{Z}S_{r,d}/S_{r,d}^{(d-1)}$ via conjugation. Moreover, this module can be naturally embedded in the linear space over the division ring (sfield) $T$ containing  $\mathbb{Z}S_{r,d}/S_{r,d}^{(d-1)}$ The dimension of this linear space is equal $n-1.$  (see for instance \cite{Romnorm}).
\end{remark}

The following statement generalize one of the statements of Theorem \ref{th:1.12}.
 \begin{theorem} (Timoshenko \cite{Timsolv}).
Let $S_{n,d}$ be a free solvable group of rank $n\geq 2$ and class $d \geq 2$. Then
  $$t(S_{n,d}) = n-1.$$
\end{theorem}

More detailed information about minimal test sets of free solvable groups can be can be obtained by the following theorem. 

\begin{theorem} (Gupta, Timoshenko \cite{Tim2}). 
Let $F_n$ be a free group of rank $n \geq 2.$ Let $R$ be a non-trivial normal subgroup of that lies in $F_n'$. Suppose that the group ring $\mathbb{Z}F_n/R$ is a domain that satisfies the Ore condition. Then for any test set $\{g_1, ..., g_{n-1}\}$ of $F_n/R'$ all elements $g_i$ belong to $R/R'$ and are linear independent over $\mathbb{Z}F_n/R$ (in other words, over the division ring of fractions of  $\mathbb{Z}F_n/R$ that exists by our assumptions).
\end{theorem}
It follows in particular that every minimal test set $T=\{g_1, ..., g_{n-1}\}$ of a free solvable group $S_{n,d}$ of rank $n$ and class $d$ lies in the last non-trivial term $S_{n,d}^{(d-1)}$ of the derived series of $S_{n,d},$  and the subgroup $H$ of Aut$S_{n,d}$ determined by $T$ coincides with $S_{n,d}^{(d-1)}$ and thus is linear. 
\begin{problem}
Let $F_n$ be a free group of rank $n\geq 3$, and $g\in F_n$ be a test element. Is the corresponding subgroup $H(g)$ linear? 
\end{problem}
Since Aut$F_2$ is linear a similar assertion for $F_2$ is obviously true. The case $F_3$ is looking as most interesting. 

If $G_{\infty}$ is a relatively free group of the countable infinite rank with the basis $Y_{\infty}=\{y_1, y_2, ..., y_i, ...\},$ the automorphism group Aut$G_{\infty}$  contains the (non-linear) subgroup $\mathbb{S}_{\infty}$ of all permutations of the countable infinite set $Y_{\infty}.$ Hence Aut$G_{\infty}$ is non-linear. Some good properties of Aut $G_{\infty}$ for the case of arbitrary metabelian or nilpotent-by-abelian variety can be found in \cite{BrRom}.

\section{On linearity of automorphism groups of relatively free algebras}

Let $\mathbb{F}$ be a field. We restrict exposition to considering the following classical algebras over $\mathbb{F}$ (the subscript $n \geq 2$ stands for the rank, i.e., the number of free generating elements): the free Lie algebra $L_n$, the free associative algebra $A_n$, the absolutely free algebra $F_n$, and the algebra $P_n$ of polynomials. It is irrelevant whether we consider these algebras with or without identity.

The automorphism group and the tame automorphism group of an algebra $C_n$ are denoted by
Aut $C_n$ and TAut$C_n$, respectively. By definition, the subgroup TAut $C_n$ is generated in Aut $C_n$ by all elementary automorphisms and all non degenerate linear changes of generators. Respectively, TAut$C_n$ is called the subgroup of {\it tame} automorphisms. By definition, the elementary automorphisms are of the form
\begin{equation}
\phi : x_i \mapsto \alpha  x_i  + f(x_1, ...,  x_{i-1}, x_{i+1}, ..., x_n), x_j \mapsto x_j  \ {\textrm{for}} \  i ≠ j, 
\end{equation}
\noindent
where $i, j = 1, ..., n, \alpha $ is a nonzero element of  $\mathbb{F}$, and $f(x_1, ..., x_{i-1}, ,... , x_{i+1}, …, x_n)$ is an element of the subalgebra generated by the generators mentioned. It is well known that if
$C_2$ coincides with $P_2$ or $A_2$  then Aut $C_2 =$ TAut $C_2$ for the case of an arbitrary field $\mathbb{F}$ (see \cite{Cz1} - \cite{M-L}). If $\mathbb{F}$ is a field of characteristic zero, then Aut $P_3 \not= $  TAut$P_3$ by  \cite{US}-\cite{SU}, and
Aut $A_3 \not= $ TAut $A_3$ by \cite{UU1}-\cite{UU2}.

An automorphism of an algebra $C_n$ is said to be {\it unitriangular} if it has the form
\begin{equation}
\varphi : x_i \mapsto x_i  + f_i(x_1, ...,  x_{i-1}),
\end{equation}
\noindent
where $i = 1, ..., n$, and $f_i(x_1, ..., x_{i−1})$ is an element of the subalgebra generated by the generators mentioned. Let TU$_n$ denote a subgroup of of TAut $C_n$  which is generated by elementary unitriangular automorphisms of the form
\begin{equation}\phi  : x_i \mapsto x_i + f(x_1 ,..., x_{i-1}), x_j \mapsto x_j \ {\textrm{for}} \  i \not= j, 
\end{equation}
\noindent
where $i, j = 1,...,n$ and $f(x_1 ,..., x_{i−1})$ is an element of the subalgebra generated by the generators mentioned. We assume that  $x_1$ is kept fixed under any unitriangular automorphism. Formally, if algebra contains constants, this element may have an image of the form $x_1 + \alpha $, where $\alpha \in \mathbb{F}$. Under such  definition, the group of unitriangular automorphisms differs inessentially from the group defined in the way indicated above, which the former contains as a subgroup. The following result will also be valid.
\begin{theorem}
 (Roman’kov, Chirkov, Shevelin \cite{RCS}). The group of tame automorphisms of the free Lie algebra $L_n$ (the free associative algebra $A_n$, the absolutely free algebra $F_n$, the algebra $P_n$ of polynomials) of rank $n \geq 4$ over a field $\mathbb{F}$ of characteristic zero admits no faithful representation by matrices over any field.
\end{theorem}

More exactly, it was  established in all these cases that the group of tame automorphisms contains a solvable subgroup of unitriangular automorphisms TU$_n$ in which the commutant of every
subgroup of finite index is not nilpotent. By a theorem of A. I. Malcev, this is impossible in matrix groups.

Similarly to TU$_n$, we define a group T$_n$ of triangular automorphisms which is generated by elementary triangular automorphisms of the form
\begin{equation}
\phi  : x_i \mapsto  \alpha x_i + f(x_1 , ..., x_{i-1}), x_j \mapsto x_j \ {\textrm{ for}} \  i \not= j, 
\end{equation}
\noindent
where $i, j = 1,...,n,$ and $\alpha $ is a nonzero element of $\mathbb{F}$. The group T$_n$ consists of all automorphisms of the form
\begin{equation}
\phi  : x_i \mapsto \alpha_i x_i + f_i(x_1 , ..., x_{i-1}), 
\end{equation}
\noindent
where $i = 1, ..., n$, $\alpha_i$ is a nontrivial element of $\mathbb{F}$, and $f_i(x_1 , ..., x_{i-1})$ are elements of the subalgebras generated by the generators specified.
\begin{theorem}
(a) (Sosnovskii \cite{S}). For $n \geq 3$ the group Aut$P_n$  over a field $\mathbb{F}$ of characteristic zero is not linear.

(b) (Bardakov, Neschadim, Sosnovskii \cite{BNS}). For $n \geq 3$ the group Aut$A_n$ over a field $\mathbb{F}$ of characteristic zero is not linear.
\end{theorem}

\begin{theorem}
(Roman’kov \cite{Rom3}). For algebras $P_n, A_n$, and $F_n$ every finitely generated subgroup $G$ of a group of triangular automorhisms admits a faithful representation by triangular matrices over $\mathbb{F}$. Consequently, the group $G$ is soluble. At the same time, every finitely generated subgroup $H$ of a group of unitriangular automorphisms admits a faithful representation by unitriangular matrices over $\mathbb{F}$. Hence the group $H$ is nilpotent.
\end{theorem}

A shorter version of this paper was published in \cite{Herald}.

\end{document}